\documentclass[11pt]{article}
\usepackage{exscale,relsize}
\usepackage{amsmath}
\usepackage{amsfonts}
\usepackage[hidelinks]{hyperref}
\usepackage[]{hyperref}
\usepackage{amssymb}
\usepackage{calc}
\usepackage{theorem}
\usepackage{pifont}      
\usepackage{array}
\usepackage{color}

\newcommand{\To}{\ensuremath{\rightrightarrows}}

\newcommand{\scal}[2]{\langle{{#1},{#2}}\rangle}
\newcommand{\Scal}[2]{\left\langle{{#1},{#2}}\right\rangle}

\newcommand{\RR}{\ensuremath{\mathbb R}}

\newcommand{\QQ}{\ensuremath{\mathbb Q}}

\newcommand{\RPP}{\ensuremath{\,\left]0,+\infty\right[}}
\newcommand{\RX}{\ensuremath{\,\left]-\infty,+\infty\right]}}
\newcommand{\RRX}{\ensuremath{\,\left[-\infty,+\infty\right]}}

\newcommand{\NN}{\ensuremath{\mathbb N}}

\newcommand{\menge}[2]{\big\{{#1} \mid {#2}\big\}}

\newcommand{\dom}{\ensuremath{\operatorname{dom}}}
\newcommand{\closu}{\ensuremath{\operatorname{cl}}}

\newcommand{\hprox}{\ensuremath{\operatorname{Prox}}}

\newcommand{\ran}{\ensuremath{\operatorname{ran}}}

\newcommand{\Fix}{\ensuremath{\operatorname{Fix}}}

\newcommand{\Id}{\ensuremath{\operatorname{Id}}}

\renewcommand{\phi}{\ensuremath{\varphi}}

\newcommand{\GF}{\ensuremath{\Gamma_{0}(X)}}

\newtheorem{theorem}{Theorem}[section]
\newtheorem{lemma}[theorem]{Lemma}
\newtheorem{fact}[theorem]{Fact}
\newtheorem{corollary}[theorem]{Corollary}
\newtheorem{proposition}[theorem]{Proposition}
\newtheorem{definition}[theorem]{Definition}

\theoremstyle{plain}{\theorembodyfont{\rmfamily}
}
\theoremstyle{plain}{\theorembodyfont{\rmfamily}
}
\theoremstyle{plain}{\theorembodyfont{\rmfamily}
}
\theoremstyle{plain}{\theorembodyfont{\rmfamily}
\newtheorem{example}[theorem]{Example}}
\theoremstyle{plain}{\theorembodyfont{\rmfamily}
\newtheorem{remark}[theorem]{Remark}}
\theoremstyle{plain}{\theorembodyfont{\rmfamily}
}

\newcommand{\pluss}{{\hskip1pt \raise1pt\vbox{\hrule width6pt \vskip1pt
\hrule width6pt}\kern-4pt{\lower1pt\hbox{\vrule height6pt \kern1pt\vrule
height6pt}}\hskip5pt}}

\newcommand{\argmin}{\mathop{\rm argmin}\limits}

\begin{document}

\title{\textsf{ Roots of the identity operator and proximal mappings: \\
(classical and phantom) cycles and gap vectors}}

\author{
Heinz H. Bauschke\thanks{Mathematics, Irving K.\ Barber School, University
of British Columbia Okanagan, Kelowna, British Columbia V1V 1V7,
Canada. E-mail: \texttt{heinz.bauschke@ubc.ca}.}
 and
Xianfu Wang\thanks{Mathematics, Irving K.\ Barber School, University
of British Columbia Okanagan, Kelowna, British Columbia V1V 1V7,
Canada. E-mail: \texttt{shawn.wang@ubc.ca}.}}

\date{February 22, 2022}

\maketitle

 \vskip 8mm

\begin{abstract} \noindent
Recently, Simons provided a lemma for a support function of a closed convex set
in a general Hilbert space and used it to prove the geometry conjecture on
cycles of projections. In this paper, we extend Simons's lemma to closed convex
functions, show its connections to Attouch--Th\'era duality, and use it to
characterize (classical and phantom) cycles and gap vectors of proximal
mappings.
\end{abstract}



\section{Introduction}

In \cite{simons}, Simons provides a new framework for studying the geometry
conjecture on cycles and gap vectors for cyclic projections; see also
\cite{salihah:d}.  His ingenious
analysis is mainly
based on
two technical lemmas: one for the support function of a nonempty closed, convex
subset, and the other for the negative displacement mapping on the null space of
an averaged operator involving the $m$th root of the identity operator.

\emph{Contributions.}
Our goal in this paper is to extend Simons's results from support functions to
proper lower semicontinuous convex functions, and
use them to
study classical and phantom
cycles and gap vectors for proximal mappings, which significantly generalize
some results in \cite{simons,salihah}. One distinguishing feature is that we can
study phantom cycles and gap
vectors of a convex function associated with an arbitrary isometry,
rather than just the right-shift operator like \cite{salihah}.

\emph{Notation and terminology.}
Notation is largely from \cite{BC2017,Simons2} to which we refer for background material on
proximal mappings, convex analysis, and monotone operator theory.
Throughout this paper, we assume that $X$ is a real Hilbert space with inner product
    $\scal{\cdot}{\cdot}\colon X\times X\to\RR$
and induced norm $\|\cdot\|=\sqrt{\scal{\cdot}{\cdot}}$.
The set of proper lower semicontinuous
convex functions from $X$ to $]-\infty,+\infty]$  is denoted by
$\Gamma_{0}(X)$. Let $f, g:X\rightarrow\RX$.
The \emph{Fenchel conjugate} of $f$ is
$$f^*:X\rightarrow \RRX: x^*\mapsto\sup_{x\in X}\big(\scal{x}{x^*}-f(x)\big).$$
The \emph{infimal convolution}
of $f$ and $g$  is $ f\Box g:X\rightarrow \RRX: x\mapsto\inf_{y\in X}(f(y)+g(x-y)),$
and it is exact at a point $x\in X$ if
$(\exists y\in X)\ (f\Box g)(x)=f(y)+g(x-y).$
The \emph{subdifferential} of $f$ is the set-valued operator
$$\partial f:X\To X: x\mapsto\menge{x^*\in X}{(\forall y\in X)\ f(y)\geq f(x)+\scal{x^*}{y-x}}.$$
For $f\in \Gamma_{0}(X)$, its \emph{proximal mapping} is defined by
$\hprox_{f}=(\Id+\partial f)^{-1}$.
We use $\closu f$ for the \emph{lower semicontinuous convex hull} of $f$.
For a set $C\subseteq X$, its \emph{indicator function} is defined by
$$\iota_{C}(x)=\begin{cases}
0, & \text{ if $x\in C$;}\\
+\infty, & \text{ if $x\not\in C$.}
\end{cases}
$$
The
closure of $C$ will be denoted by
$\overline{C}$.
When the set $C$ is convex, closed, and nonempty, then we
write $P_{C}$ for the \emph{projection operator} onto $C$ and $N_{C}=\partial \iota_{C}$
for the \emph{normal cone operator}.

An operator $N:X\rightarrow X$ is \emph{nonexpansive}
 if $(\forall x, y \in X)\ \|Nx-Ny\|\leq \|x-y\|$;
\emph{firmly nonexpansive} if $2N-\Id$ is nonexpansive; \emph{$\beta$-cocercive}
 if $\beta N$ is firmly
nonexpansive for some $\beta\in\RPP$. Prime examples of firmly nonexpansive mappings are
proximal mappings of elements of $\Gamma_{0}(X)$.
As usual, $\Fix N=\menge{x\in X}{Nx=x}$ denotes the set of fixed points of $N$.
For a set-valued operator $A: X \To X$, the sets
 $\dom A =\menge{x\in X}{Ax\neq\varnothing},
 \ran A =A(X)$,
 and $\ker A=A^{-1}(0)$
 are the \emph{domain}, \emph{range}, and \emph{kernel} of $A$ respectively.
For a linear operator $R:X\rightarrow X$, $R^*$ denotes its \emph{Hilbert adjoint}, see, e.g.,
\cite{BC2017,kreyszig}.


\emph{Background and motivation.}
Let $f\in \GF$ and $R:X\rightarrow X$ be a nonexpansive linear operator.
Every $z\in X$ satisfying
\begin{equation}\label{e:ab:cycle}
z=\hprox_f R z, \text{ equivalently, }
\end{equation}
in terms of monotone operators
\begin{equation}\label{e:sub:form}
0\in \partial f(z)+z-Rz.
\end{equation}
is called a \emph{cycle} of $f$ associated with $R$. Denote
$Z=\menge{z\in X}{z=\hprox_{f}R z}.$
The set of \emph{gap vectors} of $f$ is defined as $G=\menge{Rz-z}{z\in Z}.$
These concepts become
more meaningful and geometric, when $f$ is a decomposable sum and $R$ is the right-shift operator (see below)
 on the
product space $X^m$ with $m\in\NN=\{1, 2, \dots\}$. More precisely,
equip the product space $X^m$ with the inner product norm
$\|x\|=\sqrt{\|x_{1}\|^2+\cdots+\|x_{m}\|^2}$
for $x=(x_{1},x_{2},\ldots, x_{m})\in X^m$.
Define the right-shift operator
\begin{equation}\label{e:shifty}
R:X^m\rightarrow X^m:(x_{1},x_{2},\ldots, x_{m})\mapsto (x_{m},x_{1},\ldots, x_{m-1})
\end{equation}
and a decomposable sum of functions
\begin{equation}\label{e:dsum}
f=f_{1}\oplus\cdots\oplus f_{m}: X^m\rightarrow\RX: (x_{1},\ldots, x_{m})\mapsto \sum_{i=1}^{m}f_{i}(x_{i})
\end{equation}
where $(f_{i})_{i=1}^{m}$ in $\GF$.
A \emph{classical cycle (or proximal cycle)} of $f$ is a vector $z=(z_{1},\ldots,z_{m})\in X^{m}$ such that
\begin{subequations}
\begin{equation}\label{e:cycle}
z_{1} =\hprox_{f_{1}}z_{m},  \quad z_{2}=\hprox_{f_{2}}z_{1}, \quad z_{3}=\hprox_{f_{3}}z_{2}, \cdots,
\end{equation}
\begin{equation}\label{e:cycle:end}
z_{m-1} =\hprox_{f_{m-1}}z_{m-2},  \quad z_{m}=\hprox_{f_{m}}z_{m-1},
\end{equation}
\end{subequations}
see, e.g., \cite{salihah}.
Such a $z$ is precisely a solution to \eqref{e:ab:cycle} with $f$ and $R$ given by
\eqref{e:dsum} and \eqref{e:shifty} respectively.
In particular, for $f_{i}=\iota_{C_{i}}$ with $C_{i}$ being a
nonempty closed convex subset of $X$,
$Z$ gives the \emph{classical cycles} associated with the family of projections $P_{C_{i}}$.
See \cite{salihah:21, salihah:d, salihah,simons} for further details.

\emph{Outline.} The rest of the paper is organized as follows.
In Section~\ref{s:average} we provide some
new properties
of an averaged operator of powers of $m$th roots of the identity operator.
In Section~\ref{s:simons:lemma}
we extend Simons's lemma to lower semicontinuous convex functions
and establish its connections to Attouch--Th\'era duality.
 Section~\ref{s:classicalcg}
contains characterizations of classical cycle and gap vectors.
In the final section~\ref{s:general:cycleg} we give characterizations of
phantom cycles and gap vectors.

\section{The associated average operator: kernel and range}\label{s:average}
Let $R:X\rightarrow X$ be linear and $R^m=\Id$. Define the average operator
$$A=\frac{1}{m}\sum_{i=1}^{m}R^i,\text{
and } Y=\ker(A)=\menge{y\in X}{Ay=0}.$$
 Also define $S:X\rightarrow X$ by
$S=R-\Id$ and $Q:X\rightarrow X$ by
$Q=\frac{1}{m}\sum_{i=1}^{m-1}iR^i$, and $Q_{0}=Q|_{Y}$, the restriction of $Q$ to $Y$.
Linear operators $A$, $S$, $Q$ and subspace $Y$ are crucial in the analysis of \cite{simons}.
In this section, we show that
$A$ is in fact a projection, and that
$Y=(\Fix R)^{\perp}=\ran S$ whenever $R$ is an isometry.

We start with the following fact by Simons \cite{simons}.
We will use these properties throughout the paper.
\begin{fact}[Simons]\label{f:simons:fact} The following hold:
\begin{enumerate}
\item $S(X)\subseteq Y$, and $Q(Y)\subseteq Y$.
\item \label{i:inverse:y} $(\forall y\in Y)\ S(Qy)=y, \text{ and } Q(Sy)=y.$
\item\label{i:as} $AS=SA=0.$
\item\label{i:inverse} $SQ=QS=\Id -A$.
\item\label{i:max} $-Q_{0}-\Id/2$ is skew and so maximally monotone on $Y$.
\item \label{i:trick}
If $R$ is an isometry, then $(\forall x\in X)$ $2\scal{x}{Sx}+ \|Sx\|^2 = 0$.
\end{enumerate}
\end{fact}
\begin{example}\label{e:not:nonexpansive}
 A linear operator $R:X\rightarrow X$ satisfying $R^m=\Id$ does not imply $R$ nonexpansive.
Let $e_1,e_2,e_3,e_4$ be the canonical base of the Euclidean space $\RR^4$.
\begin{enumerate}
\item Bambaii--Chowla's matrix (1946):
Set
$$B_{1}=\begin{pmatrix}
-1 & -1& -1 & -1\\
1 &0 &0 &0\\
0& 1 &0 &0\\
0 & 0 & 1 & 0
\end{pmatrix}.
$$
Then $B_1^5=\Id$ but
$\|B_{1}e_{1}\|=\sqrt{2}>1=\|e_{1}\|$.
\item
Set
$$B_{2}=\begin{pmatrix}
1 & 1& 1 & 1\\
0 & -1 & -2 &-3\\
0& 0 &1 &3\\
0 & 0 & 0 & -1
\end{pmatrix}.
$$
Then
$B_{2}^2=\Id$ but
$\|B_{2}e_{4}\|=\sqrt{20}>1=\|e_{4}\|$.

\item Turnbull's matrix (1927):
Set
$$B_{3}=\begin{pmatrix}
-1 & 1& -1 & 1\\
-3 & 2 & -1 & 0\\
-3 & 1 &0 &0\\
-1 & 0 & 0 & 0
\end{pmatrix}.
$$
Then $B_{3}^3=\Id$ but
$\|B_{3}e_{1}\|=\sqrt{20}>1=\|e_{1}\|.$
\end{enumerate}
\end{example}
See \cite{higham} for further information on roots of matrices.
However, the following holds.
\begin{proposition}\label{p:equivalent}
 Let $R:X\rightarrow X$ be linear and $R^m=\Id$ for $m\in\NN$. Then the following are equivalent:
 \begin{enumerate}
 \item\label{i:nonex} $R$ is nonexpansive.
 \item\label{i:iso} $R$ is an isometry.
 \item \label{i:nonex1} $R^*$ is nonexpansive.
 \item \label{i:iso1} $R^*$ is an isometry.
 \end{enumerate}
\end{proposition}
\begin{proof} ``\ref{i:nonex}$\Rightarrow$\ref{i:iso}": Suppose $R$ is nonexpansive. Then $\|R\|\leq 1$.
    Using $R^m=\Id$, we obtain
$$(\forall x\in X)\ \|x\|=\|R^m x\|\leq \|R^{m-1}x\|\leq\cdots\leq \|Rx\|\leq \|x\|,$$
so $(\forall x\in X)\ \|Rx\|=\|x\|$. Hence $R$ is isometric.
``\ref{i:iso}$\Rightarrow$\ref{i:nonex}": Clear.

By the assumption, $R^*:X\rightarrow X$ is linear and $(R^*)^m=\Id$. Similar argument applying to $R^*$
yields \ref{i:nonex1}$\Leftrightarrow$\ref{i:iso1}.
Finally \ref{i:nonex}$\Leftrightarrow$\ref{i:nonex1} follows from $\|R\|=\|R^*\|$.
\end{proof}

With Example~\ref{e:not:nonexpansive} and Proposition~\ref{p:equivalent} in mind, when $R$ is an isometry
we have the following new properties of $A$ and $S$.
\begin{theorem}\label{t:a:proj} Suppose that $R$ is an isometry.
Then the following hold:
\begin{enumerate}
\item\label{i:kernel}
$\ker A=\ker A^*=(\Fix R)^{\perp}=(\Fix R^*)^{\perp}.$
\item\label{i:operator:a}
 $A=P_{\Fix R}=P_{\Fix R^*}=A^*$. In particular, $\ran A=\ran A^*=\Fix R$ is closed.
\item\label{i:operator:s} $\ran S=(\Fix R)^{\perp}=\ran S^*.$ In particular, $\ran S=\ran S^*$ is
 closed.
\end{enumerate}
\end{theorem}
\begin{proof} \ref{i:kernel}. Fact~\ref{f:simons:fact}\ref{i:as} gives
    $\ran A\subseteq\ker S=\Fix R$ and
$A^*S^*=S^*A^*=0$, so that $\ran A^*\subseteq\ker S^*=\Fix R^*$.
Because $R$ is nonexpansive and $\|R^*\|=\|R\|$, $R^*$ is nonexpansive. Then
both $\Id- R$ and $\Id-R^*$ are maximally monotone linear operators.  \cite[Proposition 3.1]{bbw07}
or \cite[Theorem 3.2(i)]{yao09} gives
$$\Fix R=\ker (\Id-R)=\ker (\Id- R^*)=\Fix R^*.$$
We have
\begin{align}
\ker A &=(\ran A^*)^{\perp} \supseteq (\Fix R^*)^{\perp}
=(\Fix R)^{\perp}.
\end{align}
To show the converse inclusion, let $y\in\ker A$. We show $y\in (\Fix R)^{\perp}$.
As $y\in\ker A$, $\sum_{i=1}^{m}R^iy=0$. For $x\in \Fix R=\Fix R^*$, we have
$\scal{R^iy}{x}=\scal{y}{(R^*)^{i}x}=\scal{y}{x}$. Then
$$0=\Scal{\sum_{i=1}^{m}R^iy}{x}=\sum_{i=1}^{m}\scal{R^iy}{x}=\sum_{i=1}^{m}\scal{y}{x}=m\scal{y}{x},
$$
i.e., $\scal{y}{x}=0$. Since this holds for every $x\in \Fix R$, we obtain
$y\in (\Fix R)^{\perp}.$
Hence $\ker A=(\Fix R)^{\perp}$.
Using a similar argument with $A$ and $A^*$ interchanged and
$R$ and $R^*$ interchanged gives
$\ker A^*=(\Fix R^*)^{\perp}=(\Fix R)^{\perp}$.

\ref{i:operator:a}. We have $\overline{\ran A}=(\ker A^*)^{\perp}=((\Fix R)^{\perp})^{\perp}=\Fix R$.
For every $x\in\Fix R$, we have $Ax=x$ and so $\Fix R\subseteq \ran A$. Thus,
$$\Fix R\subseteq \ran A\subseteq\overline{\ran A}=\Fix R,$$
which gives $\ran A=\Fix R$.
Applying this to $A^*$ yields $\ran A^*=\Fix R^*=\Fix R$.

To show $A=P_{\Fix R}$, we use $\ker A=(\Fix R)^{\perp}$ and $\ran A=\Fix R$.
For every $u\in \Fix R$, we have $Au=u$ by the definition of $A$; for every
 $v\in (\Fix R)^{\perp}$ we have $Av=0$. For each $x\in X$, by the orthogonal decomposition
 theorem, $x=u+v$ for some unique $u\in \Fix R, v\in (\Fix R)^{\perp}$. It follows that
 $$Ax=A(u+v)=Au+Av=Au=u=P_{\Fix R}u=P_{\Fix R}(u+v)=P_{\Fix R}x.$$
 Hence $A=P_{\Fix R}$ and so $A^*=P_{\Fix R}^*=P_{\Fix R}.$

 \ref{i:operator:s}. Let $y\in \ran S$. By \cite[Theorem 3.2(ii)]{yao09},
 $\overline{\ran}(\Id-R)=\overline{\ran} (\Id-R^*)$. Then
 $$y\in\overline{\ran} S=\overline{\ran} S^*=(\ker S)^{\perp}=(\Fix R)^{\perp}$$
 so $\ran S\subseteq (\Fix R)^{\perp}$.
 Conversely, let $y\in (\Fix R)^{\perp}$. As in \cite[Proposition 3.1]{salihah:21},
 setting
 $$x=\frac{1}{m}\sum_{k=0}^{m-2}(m-1-k)R^k y,$$
 we show  $y=S(-x)$.
 Indeed, using
 $A=P_{\Fix R}$, we have
 \begin{align}
 Sx & =(R-\Id) x =\frac{1}{m}(R-\Id)\sum_{k=0}^{m-2}(m-1-k)R^k y\\
 & =\left(\frac{1}{m}\sum_{k=0}^{m-1}R^k-\Id\right)y=(P_{\Fix R}-\Id) y=-P_{(\Fix R)^{\perp}}y=-y.
 \end{align}
 Hence $(\Fix R)^{\perp}\subseteq \ran S$. Altogether, $\ran S=(\Fix R)^{\perp}$.
 In view of Proposition~\ref{p:equivalent}, applying similar argument to $S^*$ yields
 $\ran S^*=(\Fix R^*)^{\perp}=(\Fix R)^{\perp}.$
\end{proof}
\begin{remark}
\begin{enumerate}
\item The referee suggested that ``$\ran S=(\Fix R)^{\perp}$" in Theorem~\ref{t:a:proj}\ref{i:operator:s}
can also be proved in the following way: By virtue of (i) it suffices to prove that
$\ran S=\ker A$. If $y\in\ran S$ then, for some $x\in X$,
$y=Sx=Rx-x$, so $Ay=ARx-Ax=Ax-Ax=0$, and $y\in\ker A$. Thus
$\ran S\subseteq\ker A$. If conversely, $y\in\ker A$ then Fact~\ref{f:simons:fact}\ref{i:inverse}
gives $y=y-Ay=SQy\in\ran S$. Thus $\ker A\subseteq\ran S$. So we have proved that
$\ran S=\ker A$, as required.
\item   The proof of Theorem~\ref{t:a:proj}\ref{i:operator:a} gives a new proof of
    the right identity in (4) of \cite[Proposition~2.4]{salihah:21}.
\end{enumerate}
\end{remark}

\begin{example} Consider the following isometric mapping $R$.
\begin{enumerate}
\item Define the right-shift operator $R:X^m\rightarrow X^m$ by
$$R(x_{1},\ldots, x_{m})=(x_{m}, x_{1},\ldots, x_{m-1}).$$
Then $R^m=\Id$, $\Fix R=\Delta$ so that $A=P_{\Delta}$ and $\ker A=\Delta^{\perp}$,
where $\Delta=\menge{(x,\ldots,x)\in X^m}{x\in X}.$

\item Define the identity operator $R:X\rightarrow X$ by $R:=\Id$.
Then $\Fix R=X$, $A=\Id$ and $\ker A=\{0\}$.

\item Define the rotator $R:\RR^2\rightarrow\RR^2$ by $R:=R_{\alpha\pi}$ where
$\alpha\in \QQ\cap\left]0,2\right[$. Let $m$ be in $\NN$ such that $m\alpha\in 2\NN$.
Then $R^m=\Id$, $\Fix R=\{0\}$, $A=0$, and $\ker A=\RR^2$.
\end{enumerate}
\end{example}

\begin{example} Without $R$ being isometric, Theorem~\ref{t:a:proj}\ref{i:operator:a} fails.
Take $B_{2}$ in Example~\ref{e:not:nonexpansive}(ii) where $m=2$ to obtain
$$A=\frac{1}{2}(B_{2}+B_{2}^2)=\begin{pmatrix}
1 & 1/2 & 1/2 & 1/2\\
0 & 0& -1 & -3/2\\
0& 0& 1& 3/2\\
0 & 0& 0 &0
\end{pmatrix}.
$$
Because $\|A e_{4}\|=\sqrt{19/4}>\|e_{4}\|$, the operator
$A$ can neither be nonexpansive nor a projection operator.
\end{example}

In view of Theorem~\ref{t:a:proj}, in the remainder of this paper, we shall assume that
$R$ is an isometry and $Y=(\Fix R)^{\perp}$.

\section{Extended Simons's lemma and Attouch--Th\'era duality}\label{s:simons:lemma}

Let $Y$, $S$, $Q$ be given as in Section~\ref{s:average}.
We call the following result the  \emph{extended Simons's lemma}. In \cite[Lemma 16]{simons},
Simons only proved this
when $f=\sigma_{C}$, a support function of a closed convex set $C\subseteq X$.
Our proof here also follows the idea of his \cite[Lemma 16]{simons}.
We also observe the uniqueness.

\begin{lemma}\label{l:simon} Let $f\in\GF$ with
$Y\cap \dom f^*\neq\varnothing$.
Then there exists a \emph{unique} pair of vectors
$(e,d)=(e_f,d_f)\in Y\times Y$ such that
$d=Se\in\dom f^*$, $e=Qd$, and
$$(\forall y\in Y)\ f^*(Se)+\scal{y-Se}{e}-f^*(y)\leq 0;$$
equivalently, $e\in\partial(f^*+\iota_Y)(Se)$.
Consequently,
$(\forall x\in X)\ f^*(Se)+\scal{Sx-Se}{e}-f^*(Sx)\leq 0.$
\end{lemma}
\begin{proof}
Set $g=f^*|_{Y}$. The assumption on $f$ implies $g\in\Gamma_{0}(Y)$ so that
$\partial g$ is maximally monotone on $Y$ by \cite[Theorem 20.25]{BC2017}.
From Fact~\ref{f:simons:fact}\ref{i:max},
$-Q_{0}-\Id/2$ is
maximally monotone on $Y$.
Since $-Q_{0}-\Id/2$ has full domain,
the operator $(-Q_{0}-\Id/2)+\partial g$, being a sum of two maximally monotone
operators, is maximally monotone on $Y$ by \cite[Corollary 25.5]{BC2017} or \cite[Theorem 1]{Rock70},
and so is $-2Q_{0}-\Id+2\partial g$.
Minty's theorem, see, e.g., \cite[Theorem 21.1]{BC2017}, implies that
there exists a \emph{unique} vector $d\in Y$ such that
$0\in -2Q_{0}d+2\partial g(d)$, i.e.,
\begin{equation}
\label{e:220111a}
Q_{0}d\in\partial g(d).
\end{equation}
Put $Q_{0}d=e$. From Fact~\ref{f:simons:fact}\ref{i:inverse:y} or
Fact~\ref{f:simons:fact}\ref{i:inverse} and the definition of $Y$ we get
$d=SQ_{0}d=Se$ and
\begin{equation}
\label{e:220111b}
(\forall y\in Y)\ g(d)+\scal{e}{y-Se}=g(d)+\scal{Q_{0}d}{y-d}\leq g(y);
\end{equation}
equivalently,
$(\forall y\in Y)\ f^*(Se)+\scal{y-Se}{e}\leq f^*(y)$
$\Leftrightarrow$
$(f^*+\iota_Y)(Se) + (f^*+\iota_Y)^*(e) \leq \scal{Se}{e}$
$\Leftrightarrow$
$e\in\partial(f^*+\iota_Y)(Se)$.
Finally, \eqref{e:220111b} is equivalent to \eqref{e:220111a} which
in turn has a \emph{unique} solution $d$ by Minty's theorem.
\end{proof}

\begin{lemma}\label{l:simon:variant}
Let $f\in\GF$ with $Y\cap \dom f^*\neq\varnothing$.
Then the vector $e=e_f\in Y$ from Lemma~\ref{l:simon} is the
\emph{unique} vector satisfying
\begin{align}
& (f^*+\iota_{Y})(Se)-\scal{Se}{e}+\closu(f\Box \iota_{Y^{\perp}})(e)\label{e:fenchel:yuong}\\
&=(f^*+\iota_{Y})(Se)-\scal{Se}{e}+(f^*+\iota_{Y})^*(e)=0.\label{e:fenchel:yuong1}
\end{align}
\end{lemma}

\begin{proof}
Lemma~\ref{l:simon} shows that $e$ is the unique vector
satisfying $e\in\partial(f^*+\iota_Y)(Se)$.
Because $Y\cap \dom f^*\neq\varnothing$, \cite[Theorem 15.1]{BC2017} implies that
$f\Box \iota_{Y^{\perp}}$ is proper, convex and bounded below by a continuous affine function,
and that
$(f^*+\iota_{Y})^*=\closu(f\Box \iota_{Y^{\perp}})$.
The result now follows from the characterization of equality in
the Fenchel--Young inequality.
\end{proof}

The extended Simons's lemma is closely related to Attouch--Th\'era duality, as we show next.
Attouch--Th\'era duality is a powerful tool in studying primal-dual solutions
of monotone inclusion problems.

\begin{fact}[Attouch--Th\'era duality \cite{AtTh}]\label{dualityat}
Let $A,B:X\To X$ be maximally monotone operators. Let
$C$ be the solution set of the primal problem:
\begin{equation}\label{theraprimal}
\text{ find $x\in X$ such that } 0\in Ax+Bx.
\end{equation}
Let $C^*$ be the solution set of the dual problem associated with the ordered pair $(A,B)$:
\begin{equation}
\label{theradual}
\text{ find $x^*\in X$ such that } 0\in A^{-1}x^*+\widetilde{B}(x^*),
\end{equation}
where
$\widetilde{B} =(-\Id)\circ B^{-1}\circ (-\Id).$
Then
\begin{enumerate}
\item $C=\menge{x\in X}{(\exists\ x^*\in C^*)\ x^*\in Ax \text{
and } -x^*\in Bx}$.
\item $C^*=\menge{x^*\in X}{(\exists\ x\in C)\ x\in A^{-1}x^*
\text{ and } -x\in \widetilde{B}(x^*)}.$
\end{enumerate}
\end{fact}
\begin{definition} We refer to \eqref{theraprimal} and \eqref{theradual} as
an Attouch--Th\'era primal-dual inclusion pair.
\end{definition}

\begin{theorem}\label{t:pd:vectors}
Let $R$ be an isometry and $Y=(\Fix R)^{\perp}$, let $f\in \GF$ with
$Y\cap \dom f^*\neq\varnothing$, and let $(e,d)\in Y\times Y$ be given by Lemma~\ref{l:simon}.
Consider the Attouch--Th\'era primal-dual inclusion problem:
\begin{align}
(P) & \quad 0\in \partial \closu(f\Box\iota_{Y^{\perp}})(x)+(\Id-R)x,\label{e:p:one}\\
(D) & \quad 0\in \partial (f^*+\iota_{Y})(y)+(\Id-R)^{-1}y. \label{e:d:one}
\end{align}
Then the following hold:
\begin{enumerate}
\item\label{i:one:sol} $(e,d)$ is a solution to the primal-dual problem \eqref{e:p:one}--\eqref{e:d:one}, i.e.,
$e$ solves $(P)$ and $d$ solves $(D)$.
Moreover, $d$ is the \emph{unique} solution of $(D)$.

\item\label{i:unique:sol}
$(e,d)$ is the \emph{unique} solution of the primal-dual problem
\begin{align}
(P') & \quad 0\in \partial \closu(f\Box\iota_{Y^{\perp}})(x)+(\Id-R)x
\;\;\text{and}\;\; x \in Y, \label{e:pu} \\
(D') & \quad 0\in \partial (f^*+\iota_{Y})(y)+(\Id-R)^{-1}y. \label{e:du}
\end{align}
More specifically,
$e$ is the \emph{unique} solution of $(P')$ and $d$ is the
\emph{unique} solution of $(D')$.
\end{enumerate}
\end{theorem}

\begin{proof}
\ref{i:one:sol}:
It is clear that \eqref{e:p:one} and \eqref{e:d:one} is an Attouch--Th\'era
primal-dual inclusion pair, because $[\partial \closu(f\Box\iota_{Y^{\perp}})]^{-1}=\partial(f^*+\iota_{Y})$ and
$\widetilde{\Id-R}=(\Id-R)^{-1}$.
We only need to verify that $(e,d)$ is a solution to the pair.
By Lemma~\ref{l:simon:variant}, $(e, d)\in Y\times Y$,
$Se=d$, and
$$(f^*+\iota_{Y})(Se)-\scal{Se}{e}+(f^*+\iota_{Y})^*(e)=0.$$
Then
$Se\in \partial(f^*+\iota_{Y})^*(e)$, that is,
$0\in (\Id-R)(e)+\partial \closu(f\Box\iota_{Y^{\perp}})(e).$
Hence $e$ solves $(P)$.

Also, $e\in \partial (f^*+\iota_{Y})(Se)=\partial (f^*+\iota_{Y})(d)$. This gives
$$0\in -e+\partial (f^*+\iota_{Y})(d).$$
Since $-e=-Qd$ and $S(Qd)=d$,
we obtain $Qd\in S^{-1}(d)$ and
$$0\in -S^{-1}(d)+\partial (f^*+\iota_{Y})(d)=(\Id-R)^{-1}(d)+\partial (f^*+\iota_{Y})(d).$$
Hence $d$ solves $(D)$.
Note that $(D)=(D')$ and we will address uniqueness in the proof of
\ref{i:unique:sol} which we tackle next.

\ref{i:unique:sol}: By \ref{i:one:sol}, \eqref{e:pu}--\eqref{e:du} has at least one solution.
It remains to prove the uniqueness.

Now the solution to $0\in \partial (f^*+\iota_{Y})(y)+(\Id-R)^{-1}y$, i.e.,
to $(D)=(D')$ is unique because
$\partial (f^*+\iota_{Y})+(\Id-R)^{-1}$ is strongly monotone: indeed, since
$(\Id-R)^{-1}=-S^{-1}$, $\dom S^{-1}=\ran S=Y$ by Theorem~\ref{t:a:proj},
and $-S^{-1}|_{Y}=-Q_{0}$
is strongly monotone on $Y$
by Fact~\ref{f:simons:fact}, we deduce that $(\Id-R)^{-1}$ is strongly monotone; or
apply the cocoercivity of
$\Id-R$, see, e.g., \cite[Fact 2.8]{salihah}.
Being a sum of a monotone operator and a strongly monotone operator,
$\partial (f^*+\iota_{Y})+(\Id-R)^{-1}$ is strongly monotone.

We have seen that $d$ is the unique solution to $(D)=(D')$.
Now let $y_1$ and $y_2$ be two solutions of $(P')$, i.e., of $(P)$
with the additional requirement that $y_1$ and $y_2$ lie in $Y$.
By Fact~\ref{dualityat}(i),
$(\Id-R)y_1 = d = (\Id-R)y_2$.
Hence
$y_1-y_2\in Y\cap \Fix R=
(\Fix R)^{\perp}\cap \Fix R=\{0\}$ and therefore
$y_{1}=y_{2}.$
\end{proof}

\begin{remark}
\begin{enumerate}
\item The referee provided a simpler proof of the uniqueness of the solution to $(D)$:
$y$ is a solution to $D$ when there exists $u\in X$ such that $(\Id-R)u=y$ and $-u
\in \partial (f^*+\iota_{Y})(y)$. Then $y=-Su$ and so $-u\in \partial (f^*+\iota_{Y})(-Su)$.
If $y'$ is also a solution to $(D)$ then there exists $u'\in X$ such that $y'=-Su'$ and
$-u'\in \partial (f^*+\iota_{Y})(-Su')$. Consequently,
$\scal{u'-u}{Su'-Su}\geq 0$. From Fact~\ref{f:simons:fact}\ref{i:trick}, $Su'=Su$, and so
$y'=y$. A similar argument shows that the solution to $(P')$ is also unique:
If $x,x'$ are solutions to $(P')$ then $Sx\in \closu(f\Box\iota_{Y^{\perp}})(x)$ and
$Sx'\in \closu(f\Box\iota_{Y^{\perp}})(x')$. Consequently,
$\scal{x'-x}{Sx'-Sx}\geq 0$. From Fact~\ref{f:simons:fact}\ref{i:trick}, $Sx'-Sx=0$, that is
to say,
$x-x'\in\ker S=\Fix R$. From Theorem~\ref{t:a:proj}\ref{i:kernel},  $x-x'\in \ker A=
(\Fix R)^{\perp}$. So $x'=x$.

\item
In \cite{salihah},
$e$ and $d$ are called ``generalized cycle'' and ``generalized gap vector'' of $f$,
respectively.
In view of Theorem~\ref{t:pd:vectors},
these vectors are
the classical cycle and gap vectors of $\closu(f\Box\iota_{Y^{\perp}})$ whenever
$Y\cap \dom f^*\neq\varnothing$.
While the solution to \eqref{e:p:one} need not be unique,
the inclusion \eqref{e:pu} always has a unique solution.
\end{enumerate}
\end{remark}

\section{Characterizations of classical cycle and gap vectors}\label{s:classicalcg}
We can use the results from Section~\ref{s:simons:lemma} to study
classical cycles and gap vectors.
While the pair $(e,d)\in Y\times Y$ given by Lemma~\ref{l:simon} always exists,
the set of classical cycle and gap vectors of $f$ might be empty;
see, e.g., \cite{salihah:d, salihah}.
We start with some elementary properties of translation-invariant functions
whose simple proofs we omit.

\begin{definition}
We say that $f:X\rightarrow\RX$ is translation-invariant with respect to
a subset $C$ of $X$ if
$f(x+c)=f(x)$ for every $x\in X$ and $c\in C$.
\end{definition}
Clearly, we have
\begin{lemma}\label{l:same}
If $f\colon X\to\RX$ is translation-invariant with respect to $C$, then
$C + \dom f\subseteq \dom f$.
\end{lemma}
\begin{lemma}\label{l:nochange} If $f\colon X\rightarrow\RX$ is translation-invariant with respect to $Y^{\perp}$, then
$f\Box\iota_{Y^{\perp}}=f.$
\end{lemma}
\begin{lemma}\label{l:infimalconv}
The following hold for every proper function $f\colon X\to\RX$:
\begin{enumerate}
\item $f\Box\iota_{Y^{\perp}}$ is translation-invariant with respect to $Y^{\perp}$.
\item
 The function
$\closu(f\Box\iota_{Y^{\perp}})$ is translation-invariant with respect to $Y^{\perp}$, namely,
$$(\forall x\in X)(\forall z\in Y^{\perp})\ \closu(f\Box\iota_{Y^{\perp}})(x+z)=\closu(f\Box\iota_{Y^{\perp}})(x).$$
\end{enumerate}
\end{lemma}
Combining Lemmas~\ref{l:nochange} and \ref{l:infimalconv} we obtain:
\begin{corollary}\label{c:double:inf}
The following hold for every function $f\colon X\to\RX$:
$$(\closu(f\Box\iota_{Y^{\perp}})) \Box \iota_{Y^{\perp}}=\closu(f\Box\iota_{Y^{\perp}}),$$
$$\closu[(\closu(f\Box\iota_{Y^{\perp}}))\Box \iota_{Y^{\perp}}]=\closu(f\Box\iota_{Y^{\perp}}).$$
\end{corollary}

Using Lemma~\ref{l:simon}, we have the following characterizations of the classical cycles of $f$.
In the proof of this theorem we shall use Fact~\ref{f:simons:fact}\ref{i:max} many times
without making explicit reference to it.

\begin{theorem}\label{t:cycle:all} Let $f\in\GF$ with $Y\cap \dom f^*\neq\varnothing$
and let $(e,d)\in Y\times Y$
    be given by Lemma~\ref{l:simon}.
Then the following statements are equivalent for every $z\in X$:
\begin{enumerate}
\item\label{i:one:n} $z=\hprox_{f} Rz$.
\item\label{i:two:n}
$f^*(Sz)+f(z)+\frac{1}{2}\|Sz\|^2=0$.
\item\label{i:three:n} $Sz=d$ and $f(z)=\closu(f\Box\iota_{Y^{\perp}})(e).$
\item\label{i:four:n} $Sz=d$ and $f(z)=\closu(f\Box\iota_{Y^{\perp}})(z).$
\end{enumerate}
\end{theorem}
\begin{proof}
\ref{i:one:n}$\Leftrightarrow$\ref{i:two:n}:
$z=\hprox_{f}Rz \Leftrightarrow Rz\in z+\partial f(z) \Leftrightarrow Sz\in \partial f(z) \Leftrightarrow$
$$f^*(Sz)+f(z)=\scal{z}{Sz}=-\frac{1}{2}\|Sz\|^2.$$

\ref{i:two:n}$\Rightarrow$\ref{i:three:n}:
By \ref{i:two},
\begin{equation}\label{e:sz:fenchel}
f^*(Sz)+f(z)+\frac{1}{2}\|Sz\|^2=0.
\end{equation}
By Lemma~\ref{l:simon},
$$f^*(Se)+\scal{Sz-Se}{e}-f^*(Sz)\leq 0.$$
Adding above two equations yields
$$f^*(Se)+f(z)+\scal{Sz-Se}{e}+\frac{1}{2}\|Sz\|^2\leq 0.$$
Since
$$f^*(Se)+f(z)\geq \scal{Se}{z},$$
by the Fenchel--Young inequality, and
$$\frac{1}{2}\|Sz\|^2=-\scal{Sz}{z},$$
we have
$$\scal{Se}{z}+\scal{Sz-Se}{e}-\scal{Sz}{z}\leq 0,$$
from which
$$-\scal{S(z-e)}{z-e}=-\scal{Sz-Se}{z-e}\leq 0.$$
Then
$\frac{1}{2}\|S(z-e)\|^2\leq 0$, so $Sz=Se=d.$
Also, by Lemma~\ref{l:simon:variant} and $\scal{Se}{e}=-\frac{1}{2}\|Se\|^2=-\frac{1}{2}\|Sz\|^2$, we obtain
\begin{equation}\label{e:se:fenchel}
f^*(Sz)+\frac{1}{2}\|Sz\|^2+\closu(f\Box \iota_{Y^{\perp}})(e)=0.
\end{equation}
Combining \eqref{e:sz:fenchel} and \eqref{e:se:fenchel} gives $f(z)=\closu(f\Box \iota_{Y^{\perp}})(e)$.

\ref{i:three:n}$\Rightarrow$\ref{i:two:n}:
Now \ref{i:three:n} ensures $Sz=d=Se$ and $\closu(f\Box \iota_{Y^{\perp}})(e)=f(z).$ Also
$\scal{Se}{e}=-\frac{1}{2}\|Se\|^2=-\frac{1}{2}\|Sz\|^2$. Then \eqref{e:fenchel:yuong}
in Lemma~\ref{l:simon:variant} gives
$$f^*(Sz)+\frac{1}{2}\|Sz\|^2+f(z) =0,$$
which is \ref{i:two}.

\ref{i:three:n}$\Leftrightarrow$\ref{i:four:n}:
Assume that
$Sz=d=Se$.
Then $z-e\in S^{-1}(0)=\Fix R$. Since
$\closu(f\Box\iota_{Y^{\perp}})$ is translation-invariant with respect to $Y^{\perp}=\Fix R$
by Lemma~\ref{l:infimalconv}(ii),
we have $\closu(f\Box\iota_{Y^{\perp}})(z)=\closu(f\Box\iota_{Y^{\perp}})(e)$.
\end{proof}

To characterize classical cycles, we have to address conditions under which
$f(z)=\closu(f\Box\iota_{Y^{\perp}})(e)$ or
$f(z)=\closu(f\Box\iota_{Y^{\perp}})(z)$. These will be investigated in the next two subsections.

\subsection{Translation-invariant functions}

\begin{lemma}\label{l:invariant}
Let $f\in\GF$ and let
$C$ be a closed linear subspace of $X$.
If $f$ is translation-invariant with respect to $C$,
then $\dom f^*\subseteq C^{\perp}$ and
$$(f^*+\iota_{C^{\perp}})^*=\closu(f\Box\iota_{C})=f\Box \iota_{C}=f.$$
\end{lemma}
\begin{proof} We can and will suppose $C^{\perp}\neq X$.
Suppose $v\not\in C^{\perp}$. Since $C$ is a subspace,
we can let $u=P_{C}v$
such that $\scal{v}{u}=\scal{v}{P_{C}v}=\|P_{C}v\|^2>0$.
Take $x_{0}\in\dom f$. Then
\begin{align*}
f^*(v) &\geq \sup_{t\in\RR}\{\scal{v}{x_{0}+tu}-f(x_{0}+tu)\}\\
&=\sup_{t\in\RR}\{\scal{v}{x_{0}}-f(x_{0})+t\scal{v}{u}\} =+\infty.
\end{align*}
Hence $\dom f^* \subseteq C^\perp$.

Next, since $\dom f^*\neq\varnothing$ and $C^{\perp}$ is a closed subspace, we have
 $\dom f^*-C^{\perp}=C^{\perp}$, so the Attouch--Brezis theorem \cite[Theorem~15.3]{BC2017} gives
 $(f^*+\iota_{C^{\perp}})^*=f\Box\iota_{C}$, which implies
 $f\Box\iota_{C}$ is lower semicontinuous, i.e., $\closu(f\Box\iota_{C})=f\Box \iota_{C}$.
 Because $C$ is a subspace and $f(x-u)=f(x)$ for $-u\in C$, we have
 $$(\forall x\in X)\ (f\Box\iota_{C})(x)=\inf_{u\in C}f(x-u)=\inf_{-u\in C}f(x-u)=\inf_{-u\in C}f(x)=f(x).$$
\end{proof}

\begin{theorem}\label{t:cycles}
Let $f\in\GF$ be translation-invariant with respect to $\Fix R$ and
such that $Y\cap \dom f^*\neq\varnothing$  where $Y=(\Fix R)^{\perp}$.
Let $d\in Y$ be given by Lemma~\ref{l:simon}.
Then the following statements are equivalent for every $z\in X$:
\begin{enumerate}
\item\label{i:one} $z=\hprox_{f} Rz$.
\item\label{i:two}
$f^*(Sz)+f(z)+\frac{1}{2}\|Sz\|^2=0$.
\item\label{i:three} $Sz=d$.
\end{enumerate}
\end{theorem}
\begin{proof}
In view of Theorem~\ref{t:cycle:all}, it is clear that
\ref{i:one}$\Leftrightarrow$\ref{i:two}$\Rightarrow$\ref{i:three}.
%
It thus suffices to show that \ref{i:three}$\Rightarrow$\ref{i:two}.
By Lemma~\ref{l:invariant}, we have
\begin{equation}\label{e:back:fin}
(f^*+\iota_{Y})^{*}=\closu(f\Box\iota_{Y^\perp})=
f\Box\iota_{Y^{\perp}}=f\Box\iota_{\Fix R}=f.
\end{equation}
Now \ref{i:three} ensures that $Sz=d=Se$, so that  $S(z-e)=0$. Then
$z-e\in\Fix R$ and $f(z)=f(e)$ by translation invariance. 
Using \eqref{e:back:fin} yields
$$(f^*+\iota_{Y})^{*}(e)=\closu(f\Box\iota_{Y^\perp})(e)=
(f\Box\iota_{Y^{\perp}})(e)=f(e)=f(z).$$
Now apply (iii)$\Rightarrow$(ii) from Theorem~\ref{t:cycle:all}.
\end{proof}

\subsection{Minimizers}

Our next result provides a sufficient condition under which the minimizers of $f$ are cycles.
More precisely, $S^{-1}(d)\cap\argmin f \subseteq \Fix (\hprox_{f}R)$ always holds.

\begin{lemma}\label{l:argmin}
Let $f\in\Gamma_0(X)$ with $Y\cap \dom f^*\neq\varnothing$
and let $(e,d)\in Y\times Y$ be given by Lemma~\ref{l:simon}.
Suppose in addition that $Sz=d$ and $z\in\argmin f$. Then
\begin{equation}\label{e:fixedpoint}
z=\hprox_{f}Rz, \text{ and }
\end{equation}
\begin{equation}\label{e:zande}
\closu(f\Box\iota_{Y^{\perp}})(e)=\closu(f\Box\iota_{Y^{\perp}})(z)=\min \closu(f\Box\iota_{Y^{\perp}})
=f(z).
\end{equation}
\end{lemma}
\begin{proof}
From Lemma~\ref{l:simon:variant},
we have
\begin{equation}\label{e:infimal}
f^*(Se)-\scal{Se}{e}+\closu(f\Box\iota_{Y^{\perp}})(e)= 0.
\end{equation}
Since
$\min f\leq \closu(f\Box\iota_{Y^{\perp}})\leq f$, we obtain
\begin{equation}\label{e:minimizer:f}
\min f=f(z)=\min \closu(f\Box\iota_{Y^{\perp}})=\closu(f\Box\iota_{Y^{\perp}})(z).
\end{equation}
Then, using $Sz=d=Se$, we obtain
\begin{align*}
0 & \leq f^*(Sz)+\tfrac{1}{2}\|Sz\|^2+f(z) =f^*(Se)+\tfrac{1}{2}\|Se\|^2 +f(z)\\
 & = f^*(Se)-\scal{Se}{e} +f(z)=-\closu(f\Box\iota_{Y^{\perp}})(e)+f(z)\\
 &= -\closu(f\Box\iota_{Y^{\perp}})(e)+\min \closu(f\Box\iota_{Y^{\perp}})\leq 0,
 \end{align*}
 which in turn implies
 \begin{equation}\label{e:backcycle}
 0=f^*(Sz)+\tfrac{1}{2}\|Sz\|^2+f(z) \;\;\text{and}\;\;
 \closu(f\Box\iota_{Y^{\perp}})(e)=f(z).
 \end{equation}
 Hence, \eqref{e:fixedpoint} follows from Theorem~\ref{t:cycle:all},
 and \eqref{e:zande} follows from \eqref{e:minimizer:f} and \eqref{e:backcycle}.
\end{proof}

\begin{theorem}\label{t:cycles:minimum}
Let $f\in\Gamma_0(X)$ with $Y\cap \dom f^*\neq\varnothing$ and
let $d\in Y$ be given by Lemma~\ref{l:simon}.
Then the following statements are equivalent for every $z\in \argmin f$:
\begin{enumerate}
\item\label{i:one} $z=\hprox_{f} Rz$.
\item\label{i:two}
$f^*(Sz)+\tfrac{1}{2}\|Sz\|^2+f(z)=0$.
\item\label{i:three} $Sz=d$.
\end{enumerate}
\end{theorem}
\begin{proof}
Combine Lemma~\ref{l:argmin} and Theorem~\ref{t:cycle:all}.
\end{proof}

Immediately we obtain the following result of
Simons \cite[Theorem 7]{simons}.
\begin{corollary}\label{t:cycles:indic}
Let $C$ be a nonempty closed convex subset of $X$.
Let $d\in Y$ be given by Lemma~\ref{l:simon} with $f=\iota_C$.
Then the following statements are equivalent for every $z\in C$:
\begin{enumerate}
\item\label{i:one} $z=P_{C} Rz$.
\item\label{i:two} $\sigma_{C}(Sz)+\frac{1}{2}\|Sz\|^2=0$.
\item\label{i:three} $Sz=d$.
\end{enumerate}
\end{corollary}
\begin{proof} Note that
    $f^*=\sigma_C$ and $0\in Y\cap \dom \sigma_C$.
Noting that $C=\argmin f$, we observe that the result is clear from
Theorem~\ref{t:cycles:minimum}.
\end{proof}

\section{Phantom cycles and phantom gap vectors}\label{s:general:cycleg}

The next result makes it clear that the classical cycles and gap vector of a function $f$
are closely related to those of $\closu(f\Box \iota_{Y^\bot})$ and to
which we refer as \emph{phantom cycles} and \emph{phantom gap vector}.


%

\begin{theorem}\label{t:gencycle}
Let $f\in \GF$
with $Y\cap \dom f^*\neq\varnothing$ and
let $(e,d) = (e_f,d_f)$ be given by Lemma~\ref{l:simon}.
Then the following hold:
\begin{enumerate}
\item
The set $Z$ of \emph{phantom cycles} of $f$, which are defined to be
the set of classical cycles of the function $\closu(f\Box \iota_{Y^\bot})$, i.e.,
$
Z=\menge{z\in X}{z=\hprox_{\closu(f\Box\iota_{Y^{\perp}})}(Rz)}$, is always nonempty
and $Z=e+Y^{\perp}.$ Consequently, $Z$ contains infinitely many elements whenever
$Y^{\perp}=\Fix R\neq \{0\}.$
\item
The \emph{phantom gap vector} of $f$,
i.e., the gap vector $d_{\closu(f\Box \iota_{Y^\bot})}$, is equal to $d=Sz\in Y$
for every $z\in Z$; moreover,
$e_{\closu(f\Box \iota_{Y^\bot})}=e$.
\end{enumerate}
\end{theorem}
\begin{proof}
Consider the function $\closu(f\Box \iota_{Y^\bot})$.
This function belongs to $\GF$, and its Fenchel conjugate is
$f^*+\iota_Y$.
Moreover, $\varnothing\neq Y \cap \dom f^* = \dom \iota_Y \cap \dom f^*
= \dom(f^*+\iota_Y) = Y\cap \dom[\closu(f\Box \iota_{Y^\bot})]^*$.
We thus may and do apply Lemma~\ref{l:simon} with $f$ replaced by $\closu(f\Box \iota_{Y^\bot})$
to obtain the two ``phantom'' vectors
$(e',d')=(e_{\closu(f\Box \iota_{Y^\bot})},d_{\closu(f\Box \iota_{Y^\bot})})
\in Y\times Y$.
Applying Lemma~\ref{l:simon:variant} to $\closu(f\Box\iota_{Y^{\perp}})$, we
learn that
\begin{equation}
\label{e:220111c}
\Big([\closu(f\Box\iota_{Y^{\perp}})]^*+\iota_{Y}\Big)(Se')
+ \Big(\closu\big([\closu(f\Box\iota_{Y^{\perp}})]\Box\iota_{Y^\perp}\big)\Big)(e')
-\scal{Se'}{e'} = 0.
\end{equation}
In view of
$[\closu(f\Box\iota_{Y^{\perp}})]^* = f^*+\iota_{Y}$ and
Corollary~\ref{c:double:inf},
we see that \eqref{e:220111c} simplifies to
\begin{equation}
\label{e:220111d}
f^*(Se') +
\big(\closu(f\Box\iota_{Y^{\perp}})\big)(e')
-\scal{Se'}{e'} = 0.
\end{equation}
Good news! Because $e'\in Y$, we deduce from
the uniqueness assertion of Lemma~\ref{l:simon:variant} that $e'=e$.
It follows that $d'=Se'=Se=d$.
Theorem~\ref{t:cycles} applied to $\closu(f\Box\iota_{Y^{\perp}})$ gives
$Z=S^{-1}d'=S^{-1}d=e+\ker S= e+\Fix R=e+Y^{\perp}$.
Finally, $(\forall z\in Z)$
$Sz \in S(Z)=\{d'\}=\{d\}$ and we are done.
\end{proof}

\begin{corollary} Let $C$ be a nonempty closed convex subset of $X$.
Then the following hold:
\begin{enumerate}
\item
The set of phantom cycles of $\iota_{C}$, i.e.,
$Z=\menge{z\in X}{z=\hprox_{\iota_{\overline{C+Y^{\perp}}}}(Rz)}$, is always nonempty
and $Z=e+Y^{\perp}$,
where $e\in \overline{C+Y^{\perp}}\cap Y$ and
$\sigma_{C}(Se)-\scal{Se}{e}=0$.
Consequently, $Z$ contains infinitely many elements as long as
$Y^{\perp}=\Fix R\neq \{0\}.$
\item
 The unique phantom gap vector of $\iota_{C}$ is $d=Sz\in Y$
for every $z\in Z$.
\end{enumerate}
\end{corollary}
\begin{proof}
Set $f=\iota_C$ and note that $f^*=\sigma_C$ and
$0\in Y\cap \dom f^*$.
Moreover,
$\closu(f\Box\iota_{Y^{\perp}}) =
\closu(\iota_{C}\Box\iota_{Y^{\perp}})=\closu{(\iota_{C+Y^{\perp}}})=
\iota_{\overline{C+Y^{\perp}}}.$
The conclusion thus follows from Theorem~\ref{t:gencycle} and
Lemma~\ref{l:simon:variant}.
\end{proof}

\begin{remark} Theorem~\ref{t:gencycle} generalizes \cite[Theorem 4.9]{salihah},
where only $R:X^m\rightarrow X^m$ given by $R(x_{1},\ldots, x_{m})=(x_{m}, x_{1},\ldots, x_{m-1})$
is considered.
\end{remark}

%
%
%
%

\section*{Acknowledgments}
The authors
would like to thank the referee for
careful reading of the manuscript and valuable suggestions, which improved the
exposition considerably.
HHB and XW were supported by NSERC Discovery grants.

\end{document}